\RequirePackage{ifpdf}
\ifpdf % We are running pdfTeX in pdf mode
\documentclass[pdftex]{sigma}
\else
\documentclass{sigma}
\fi

\usepackage{bbm}

\renewcommand{\l}{\lambda}

\def\mf{\mathfrak}
\newcommand{\al}{\alpha}
\newcommand{\be}{\beta}
\newcommand{\ga}{\gamma}
\newcommand{\de}{\delta}
\newcommand{\ep}{\varepsilon}
\newcommand{\et}{\eta}
\newcommand{\tf}{\vartheta}
\newcommand{\om}{\omega}

\begin{document}
\allowdisplaybreaks

\renewcommand{\PaperNumber}{004}

\FirstPageHeading

\renewcommand{\thefootnote}{$\star$}

\ShortArticleName{Construction of the Bethe State for the $E_{\tau,\eta}(so_3)$
Elliptic Quantum Group}

\ArticleName{Construction of the Bethe State\\ for the $\boldsymbol{E_{\tau,\eta}(so_3)}$ Elliptic 
Quantum Group\footnote{This paper is a contribution 
to the Proceedings of the Workshop on 
Geometric Aspects of Integ\-rable Systems
 (July 17--19, 2006, University of Coimbra, Portugal).
The full collection is available at 
\href{http://www.emis.de/journals/SIGMA/Coimbra2006.html}{http://www.emis.de/journals/SIGMA/Coimbra2006.html}}}

\Author{Nenad MANOJLOVI\'C and Zolt\'an NAGY}

\AuthorNameForHeading{N. Manojlovi\'c and Z. Nagy}

\Address{Departamento de Matem\'atica, FCT, Universidade do Algarve,\\ Campus de Gambelas, 8005-139 Faro, Portugal}
\Email{\href{mailto:nmanoj@ualg.pt}{nmanoj@ualg.pt}, \href{mailto:znagy@ualg.pt}{znagy@ualg.pt}}

\ArticleDates{Received October 31, 2006, in f\/inal form
December 28, 2006; Published online January 05, 2007}

\Abstract{Elliptic quantum groups can be associated to solutions of the star-triangle relation of statistical mechanics.
In this paper, we consider the particular case of the $E_{\tau,\eta}(so_3)$ elliptic quantum group.
In the context of algebraic Bethe ansatz, we construct the corresponding Bethe creation operator for the
transfer matrix def\/ined in an arbitrary representation of $E_{\tau,\eta}(so_3)$.}


\Classification{82B23; 81R12; 81R50}

\section{Introduction}

In this article, we explain the f\/irst step towards the application of \emph{algebraic} Bethe ansatz to
the elliptic (or dynamical) quantum group $E_{\tau,\eta}(so_3)$. The elliptic quantum group is
the algebraic structure associated to elliptic solutions of the star-triangle relation
which appears in interaction-round-a-face models in statistical mechanics. As it was shown
by Felder \cite{Fe}, this structure is also related to Knizhnik--Zamolodchikov--Bernard equation of conformal
f\/ield theory on tori. In fact, to each solution of the (see \cite{Ji}) star-triangle relation
a dynamical $R$-matrix can be associated. This $R$-matrix, in turn, will def\/ine
an algebra similar to quantum groups appearing in the quantum inverse scattering method (QISM).
Despite all the dif\/ferences,
this new structure preserves a prominent feature of quantum groups: a tensor product of representations
can be def\/ined.

The adjective dynamical refers to the fact that the $R$-matrix appearing in these
structures contains a parameter which in the classical limit will be interpreted as the position coordinate
on the phase space of a classical system and the resulting classical $r$-matrix will depend on it. In
the quantum setting, apart from the appearance of this extra parameter the Yang--Baxter equation (YBE) is
also deformed. At the technical level, the main dif\/ference between usual quantum groups and the one
we are about to describe lies not so much in the elliptic nature of the appearing functions as rather in
the introduction of the extra ``dynamical'' parameter and the corresponding deformation the YBE.

In QISM, the physically interesting quantity is the transfer matrix. The Hamiltonian of the
model and other observables are derived from it. The knowledge of its spectrum is thus essential.
Dif\/ferent kinds of methods under the federating name of Bethe
ansatz have been developed to calculate the eigenvalues of the transfer matrix \cite{Fa,Ko,Ku}.
The question whether the algebraic Bethe ansatz (ABA) technique can be applied to transfer matrices
appearing in the context of dynamical quantum groups has received an af\/f\/irmative answer from Felder and
Varchenko \cite{Fe2,FeVa}.
They showed how to implement ABA for the elliptic quantum group $E_{\tau,\eta}(sl_2)$, they also showed its
applications to IRF models and Lam\'e equation. Later, for the $E_{\tau,\eta}(sl_n)$ elliptic
quantum group the nested Bethe ansatz method was used \cite{Es,Sa}
and a relation to Ruijsenaars--Schneider~\cite{Sa} and quantum Calogero--Moser Hamiltonians was established \cite{ABB}.

In the f\/irst section we introduce the basic def\/initions of dynamical $R$-matrix, Yang--Baxter equation, representations,
operator algebra and commuting transfer matrices.
We def\/ine ele\-ments~$\Phi_n$ in the operator algebra which have the necessary symmetry properties to be the
creation operators of the corresponding Bethe states. As it turns out, the creation operators are not
simple functions of the Lax matrix entries, unlike in \cite{FeVa}, but they are complicated
polynomials of three generators $A_1(u)$, $B_1(u)$, $B_2(u)$ in the elliptic operator algebra. We give the
recurrence relation which  def\/ines the creation operators. Moreover, we give explicit formulas of
the eigenvalues and Bethe equations for $n = 1,2$. This strongly suggests that for higher $n$ these are correct choice
of creation operators.
 Derivation of the eigenvalues and the corresponding Bethe equations for general $n$  (from the
usual cancelation of the unwanted terms) will be published elsewhere.

\section[Representations of $E_{\tau,\eta}(so_3)$ and transfer matrices]{Representations 
of $\boldsymbol{E_{\tau,\eta}(so_3)}$ and transfer matrices}

\subsection[Definitions]{Def\/initions}
Let us f\/irst recall the basic def\/initions which will enter our construction. First, we f\/ix two complex numbers
$\tau$, $\eta$ such that ${\rm Im}(\tau) > 0$.
The central object in this paper is the $R$-matrix $R(q,u)$ which depends on two arguments $q,u \in \mathbb{C}$:
the f\/irst one is referred to
as the dynamical parameter, the second one is called the spectral parameter. The elements of the $R$-matrix
are written in terms of Jacobi's theta function:
\[
\tf(u)=-\sum_{j\in \mathbb{Z}}\exp \left(\pi i\left(j+\frac{1}{2}\right)^2\tau+2\pi i \left(j+\frac{1}{2}\right)
\left(u+\frac{1}{2}\right)\right).
\]
This function has two essential properties. It is quasiperiodic:
\[
\tf(u+1)=-\tf(u), \qquad \tf(u+\tau)=-e^{-i \tau-2 i u}\tf(u)
\]
and it verif\/ies the identity:
\begin{gather*}
\tf(u+x)\tf(u-x)\tf(v+y)\tf(v-y)=\tf(u+y)\tf(u-y)\tf(v+x)\tf(v-x)\\
\qquad{}+\tf(u+v)\tf(u-v)\tf(x+y)\tf(x-y).
\end{gather*}

The entries of the $R$-matrix are written in terms of the following functions
\begin{gather*}
g(u)  =  \frac{\tf(u-\et) \tf(u-2\et)}{\tf(\et) \tf(2 \et)},\\
\al(q_1,q_2,u)  =  \frac{\tf(\et-u)\tf(q_{12}-u)}{\tf(\et)\tf(q_{12})},\\
\be(q_1,q_2,u) =  \frac{\tf(\et-u)\tf(u)\tf(q_{12}-2\et)}{\tf(-2\et)\tf(\et)\tf(q_{12})},\\
\ep(q,u)= \frac{\tf(\et+u)\tf(2\et-u)}{\tf(\et)\tf(2\et)}-\frac{\tf(u)\tf(\et-u)}{\tf(\et)\tf(2\et)}
\left( \frac{\tf(q+\et)\tf(q-2\et)}{\tf(q-\et)\tf(q)}+\frac{\tf(q-\et)\tf(q+2\et)}{\tf(q+\et)\tf(q)}\right),\\
\ga(q_1,q_2,u)=  \frac{\tf(u)\tf(q_1+q_2-\et-u)\tf(q_1-2\et)\tf(q_2+\et)}{\tf(\et)\tf(q_1+q_2-2\et)\tf(q_1+\et)\tf(q_2)},\\
\de(q,u)  =  \frac{\tf(u-q)\tf(u-q+\et)}{\tf(q)\tf(q-\et)}.
\end{gather*}

The $R$-matrix itself will act on the tensor product $V \otimes V$ where $V$ is a three-dimensional
complex vector space with the standard basis $\{e_1,e_2,e_3\}$. The matrix units $E_{ij}$ are def\/ined
in the usual way: $E_{ij}e_k=\delta_{jk}e_i$. We will also need the following diagonal matrix 
later on $h=E_{11}-E_{33}$.

Now we are ready to write the explicit form of the $R$-matrix. The matrix is obtained via a~gauge transformation
from the solution
of the star-triangle relation which is associated to the vector representation of $B_1$ \cite{Ji}. According to
a remark in \cite{Ji}, that solution can also be derived as a symmetric tensor product (i.e. fusion) of the $A_1$
solution
\begin{gather}
R(q,u)=g(u)E_{11}\otimes E_{11}+g(u)E_{33}\otimes E_{33}+\ep(q,u)E_{22}\otimes E_{22}\nonumber\\
\phantom{R(q,u)=}{}+\al(\eta,q,u)E_{12}\otimes E_{21}+\al(q,\eta,u)E_{21}\otimes E_{12}
+\al(-q,\eta,u)E_{23}\otimes E_{32}\nonumber\\
\phantom{R(q,u)=}{}+ \al(\eta,-q,u)E_{32}\otimes E_{23}\nonumber\\
\phantom{R(q,u)=}{}+ \be(\eta,q,u) E_{11}\otimes E_{22}+\be(q,\eta,u) E_{22}\otimes E_{11}
+\be(-q,\eta,u) E_{22}\otimes E_{33}\nonumber\\
\phantom{R(q,u)=}{}+ \be(\eta,-q,u)E_{33}\otimes E_{22}\nonumber\\
\phantom{R(q,u)=}{}+\ga(-q,q,u)E_{11}\otimes E_{33}+\ga(-q,\eta,u)E_{12} \otimes E_{32}
- \ga(\eta,q,u) E_{21} \otimes E_{23}\nonumber\\
\phantom{R(q,u)=}{}+ \ga(q,-q,u) E_{33} \otimes E_{11}+ \ga(q,\eta,u) E_{32} \otimes E_{12}
- \ga(\eta,-q,u) E_{23} \otimes E_{21}\nonumber\\
\phantom{R(q,u)=}{}+ \de(q,u) E_{31}\otimes E_{13}+\de(-q,u) E_{13} \otimes E_{31}.\label{Rmat}
\end{gather}
This $R$-matrix also enjoys the unitarity property:
\[
R_{12}(q,u)R_{21}(q,-u)=g(u)g(-u)\mathbbm{1}
\]
and it is of zero weight:
\[
\left[h \otimes \mathbbm{1}+\mathbbm{1} \otimes h,R_{12}(q,u)\right]=0 \qquad  (h \in \mf{h}).
\]

The $R$-matrix also obeys the dynamical quantum Yang--Baxter equation (DYBE) in ${\rm End}(V \otimes V \otimes V)$:
\begin{gather*}
R_{12}(q-2\eta h_3,u_{12})  R_{13}(q,u_1) R_{23}(q-2\eta h_1,u_2)\\
\qquad{}= R_{23}(q,u_2)R_{13}(q-2\eta h_2,u_1)
R_{12}(q,u_{12}),
\end{gather*}
where the ``dynamical shift'' notation has the usual meaning:
\begin{gather}\label{shift}
R_{12}(q-2\eta h_3,u) \cdot v_1\otimes v_2 \otimes v_3 = \left(R_{12}(q-2\eta \lambda,u) v_1\otimes v_2\right)
\otimes v_3,
\end{gather}
whenever $h v_3= \lambda v_3$. Shifts on other spaces are def\/ined in an analogous manner.
Notice that the notion and notation of ``dynamical shift'' can be extended
to dif\/ferent situations as well, even if the appearing (possibly dif\/ferent) vector spaces $V_i$
are not 3-dimensional. For this, one only needs to
verify two conditions: that an action of $h$ is def\/ined on each $V_i$, and that each $V_i$ is a direct sum of
the weight subspaces $V_i[\lambda]$ def\/ined by that action of $h$. It is easy to see then that equation (\ref{shift})
makes sense.
Furthermore, along these lines the notion
of dynamical quantum Yang--Baxter equation and of the corresponding algebraic structures can be extended to the case
where $h$ is replaced by a higher rank Abelian Lie algebra $\mathfrak{h}$. However, in this paper
we only deal with the a special rank-one case, so from now on $\mathfrak{h}=\mathbb{C}h$. It will be clear
from the context how to generalize the relevant notions to the higher rank case.

Let us also describe a more intuitive way of looking at this shift. Def\/ine f\/irst
the shift operator acting on functions of the
dynamical parameter:
\[
\exp(2\eta \partial_q) f(q)=f(q+2\eta) \exp(2\eta \partial_q).
\]
Then equation (\ref{shift}) can also be written in the following form:
\[
R_{12}(q-2\eta h_3,u) = \exp(-2\eta h_3 \partial_q) R_{12}(q,u) \exp(2\eta h_3 \partial_q)
\]
in the sequel we will use whichever def\/inition is the f\/ittest for the particular point in our calculation.

\subsection{Representation, operator algebra}

Now we describe the notion of representation of (or module over) $E_{\tau,\eta}(so_3)$. It
 is a pair $(\mathcal{L}(q,u),W)$ where $W$ is a diagonalizable $\mf{h}$-module, that is, $W$ is a direct sum of
 the weight subspaces
$W=\oplus_{\l \in \mathbb{C}}W[\l]$ and $\mathcal{L}(q,u)$ is an operator in $\mathrm{End}(V \otimes W)$ obeying:
\begin{gather}
R_{12}(q-2\eta h_3,u_{12})  \mathcal{L}_{13}(q,u_1) \mathcal{L}_{23}(q-2\eta h_1,u_2)\nonumber\\
\qquad {}= \mathcal{L}_{23}(q,u_2)\mathcal{L}_{13}(q-2\eta h_2,u_1)
R_{12}(q,u_{12}).\label{RLL}
\end{gather}

$\mathcal{L}(q,u)$ is also of zero weight
\[
\left[h_V \otimes \mathbbm{1}+\mathbbm{1} \otimes h_W , \mathcal{L}_{V,W}(q,u)\right]=0 \qquad (h \in \mf{h}),
\]
where the subscripts remind the careful reader that in this formula $h$ might act in a dif\/ferent way on spaces
$W$ and $V$.

An example is given forthwith by $W=V$ and $\mathcal{L}(q,u)=R(q,u-z)$ which is called the fundamental
representation with evaluation point $z$.
A tensor product of representations can also be def\/ined which corresponds to the existence of a coproduct-like
structure at the abstract algebraic level. Let $(\mathcal{L}(q,u),X)$ and $(\mathcal{L}'(q,u),Y)$
be two $E_{\tau,\eta}(so_3)$
modules, then $\mathcal{L}_{1X}(q-2\eta h_Y,u)\mathcal{L}_{1Y}(q,u)$, $X\otimes Y$ is a
representation of $E_{\tau,\eta}(so_3)$ on
$X \otimes Y$ endowed, of course, with the tensor product $\mf{h}$-module structure.

The operator $\mathcal{L}$ is reminescent of the quantum Lax matrix in the FRT formulation
of the quantum inverse scattering
method, although it obeys a dif\/ferent exchange relation,
therefore we will also call it a Lax matrix. This allows us to view
the $\mathcal{L}$ as a matrix with operator-valued entries.

Inspired by that interpretation, for any module over $E_{\tau,\eta}(so_3)$ we def\/ine the corresponding
\textit{operator algebra}.
Let us take an arbitrary representation $\mathcal{L}(q,u) \in \mathrm{End}(V \otimes W)$.
The elements of the operator algebra corresponding to this representation will act on the space $\mathrm{Fun}(W)$ of
meromorphic functions of $q$ with values in $W$. Namely let $L \in \mathrm{End}(V \otimes \mathrm{Fun}(W))$
be the operator def\/ined as
\[
L(u)=\left( \begin{array}{ccc}
A_1(u)& B_1(u)& B_2(u)\\
C_1(u) & A_2(u) & B_3(u)\\
C_2(u) & C_3(u) &A_3(u)
 \end{array}\right)=\mathcal{L}(q,u)e^{-2\eta h \partial_q}.
\]
We can view it as a matrix with entries in $\mathrm{End}(\mathrm{Fun}(W))$.
It follows from equation (\ref{RLL}) that $\tilde{L}$
verif\/ies:
\begin{gather}\label{RLLti}
R_{12}(q-2\eta h,u_{12})  L_{1W}(q,u_1) L_{2W}(q,u_2)= L_{2W}(q,u_2)L_{1W}(q,u_1) 
\tilde{R}_{12}(q,u_{12})
\end{gather}
with $\tilde{R}_{12}(q,u):= \exp(2\eta(h_1+h_2)\partial_q)R_{12}(q,u)\exp(-2\eta(h_1+h_2)\partial_q)$.

The zero weight condition on $L$ yields the relations:
\begin{gather*}
\left[h,A_i\right]=0 ,\qquad  \left[h,B_j\right]=-B_j \quad (j=1,3), \qquad \left[h,B_2\right]=-2B_2,\\
\left[h,C_j\right]=C_j \quad (j=1,3), \qquad \left[h,C_2\right]=2C_2,
\end{gather*}
so $B_i$'s act as lowering and $C_i$'s as raising operators.

And f\/inally the following theorem shows how to associate a family of commuting quantities to a representation
of the elliptic quantum group.

\begin{theorem}
Let $W$ be a representation of $E_{\tau,\eta}(so_3)$. Then the transfer matrix def\/ined by $t(u)={\rm Tr}\, \tilde{L}(u) \in
\mathrm{End}(\mathrm{Fun}(W))$
preserves the subspace $\mathrm{Fun}(W)[0]$ of functions with values in the zero weight subspace of $W$.
When restricted to this subspace, they commute at different values of the spectral parameter:
\[
\left[t(u),t(v)\right]=0.
\]
\end{theorem}

\begin{proof}
The proof is analogous to references \cite{FeVa3,ABB}.
\end{proof}

The eigenvalues of the transfer matrix can be found by means of the algebraic Bethe ansatz.
In the next section we develop the f\/irst steps in this direction.
\section{Construction of the Bethe state}
\subsection{The framework of the algebraic Bethe ansatz}

The above theorem tells us how to associate the transfer matrix to an arbitrary representation of the
dynamical quantum group. Our aim is to determine the spectrum of such a transfer matrix in the usual
sense of the Bethe ansatz techniques, i.e.\ to write the Bethe equations f\/ixing the eigenvalues.

In order for the algebraic Bethe ansatz to work, this representation must
be a highest weight representation, that is possess a highest weight vector $|0\rangle$
(also called pseudovacuum) which is annihilated by the raising operators and is an eigenvector
of the diagonal elements of the quantum Lax matrix
\[
C_i(u)|0\rangle=0, \qquad A_i(u)|0\rangle=a_i(q,u)|0\rangle, \qquad i=1,2,3.
\]
Actually, any vector of the form $|\Omega \rangle = f(q)|0\rangle$ is also a highest weight vector of the representation
in question. This freedom in choosing the highest weight vector will prove essential in the sequel, so we
do not f\/ix the arbitrary function $f(q)$ for the moment. The preceding relations are modif\/ied as follows:
\begin{gather*}
C_i(u)|\Omega\rangle=0, \quad i=1,2,3, \qquad A_1(u)|\Omega\rangle=a_1(q,u)\frac{f(q-2\eta)}{f(q)}|\Omega\rangle,\\
A_2(u)|\Omega\rangle=a_2(q,u)|\Omega\rangle, \qquad A_3(u)|\Omega\rangle=a_3(q,u)\frac{f(q+2\eta)}{f(q)}|\Omega\rangle.
\end{gather*}

The representations obtained by tensorising the fundamental vector representation possesses this highest weight
vector and its transfer matrix is the transfer matrix of an IRF (interaction-round-a-face) model
with Boltzmann weights derived from the dynamical $R$-matrix (\ref{Rmat}). In this case we also have the property
that $a_1(q,u)$ does not depend on $q$, this is what we will assume in the sequel. Other representation are related
to Lam\'e equation or Ruijsenaars--Schneider Hamiltonians (see \cite{FeVa} for the $E_{\tau,\eta}(sl_2)$ case).
This is expected to happen in the $E_{\tau,\eta}(so_3)$ case, too, and we hope to report on progress in representations
and related models soon.

Once the pseudovacuum is f\/ixed, one looks for eigenvalues in the form:
\[
\Phi_n(u_1,\ldots,u_n)|\Omega\rangle
\]
under some simple (symmetry) assumptions on the lowering operator $\Phi_n$.
In the XXZ model, or for $E_{\tau,\eta}(sl_2)$, $\Phi_n$ is a simple product of the only lowering operator
$B(u)$. We will explain later, in analogy with the Izergin--Korepin model,
why $\Phi_n$ is not that simple in the $E_{\tau,\eta}(so_3)$ case.

The main result of this paper is the construction of $\Phi_n$ 
(the Bethe state) for the $E_{\tau,\eta}(so_3)$ dynamical
quantum group under simple assumptions.

Finally,  one calculates the action of the transfer matrix on the Bethe state. This will yield
3 kinds of terms. The f\/irst part (usually called wanted terms in the literature)
will tell us the eigenvalue of the transfer matrix, the second part (called unwanted terms) must be annihilated
by a careful choice of the spectral parameters $u_i$ in $\Phi_n(u_1,\ldots,u_n)$; the vanishing of these unwanted
terms is ensured if the $u_i$ are solutions to the so called Bethe equation. The third part contains terms
ending with a raising operator acting on the pseudovacuum and thus vanishes. We hope to report soon on
the form of the Bethe equations and eigenvalues, too.

Right now, we propose to develop step 2 and write the recurrence relation def\/ining $\Phi_n$.
We thus assume that a representation with highest weight vector pseudovacuum already exists.

\subsection{The creation operators}
We explicitly write the commutation relations coming from the $RLL$ relations (\ref{RLLti})
which will be used in the construction of the Bethe state
\begin{gather}
B_1(u_1)B_1(u_2)=\omega_{21}\left(B_1(u_2)B_1(u_1)-\frac{1}{y_{21}(q)}B_2(u_2)A_1(u_1)\right)+
\frac{1}{y_{12}(q)}B_2(u_1)A_1(u_2), \label{crB1B1}  \\
A_1(u_1)B_1(u_2)=z_{21}(q)B_1(u_2)A_1(u_1)-\frac{\al_{21}(\eta,q)}{\be_{21}(q,\eta)}B_1(u_1)A_1(u_2), 
\label{crB1B1+} \\
A_1(u_1)B_2(u_2)=\frac{1}{\ga_{21}(q,-q)}\left( g_{21}B_2(u_2)A_1(u_2)+\ga_{21}(\eta,-q)B_1(u_1)B_1(u_2) \nonumber
\right.\\
\left.\phantom{A_1(u_1)B_2(u_2)=}{} -\de_{21}(-q)B_2(u_1)A_1(u_1)\right), \label{crB1B1++}\\
 B_1(u_2)B_2(u_1)=\frac{1}{g_{21}}\left( \be_{21}(-q,\eta)B_2(u_1)B_1(u_2)+\al_{21}(\eta,-q)B_1(u_1)B_2(u_2)\right),
 \label{crB1B1+++}
 \\
 B_2(u_2)B_1(u_1)=\frac{1}{g_{21}}\left( \be_{21}(\eta,-q)B_1(u_1)B_2(u_2)
 +\al_{21}(-q,\eta)B_2(u_1)B_1(u_2)\right),\label{crB2B1}
\end{gather}
where
\begin{gather*}
\omega(q,u)=\frac{\ep(q,-u) \ga(q,-q,-u)+\ga(q,\eta,-u)\ga(\eta,-q,-u)}{g(-u)\ga(q,-q,-u)},\\
y(q,u)=\frac{\ga(q,-q,u)}{\ga(q,\eta,u)},\qquad
z(q,u)=\frac{g(u)}{\be(q,\eta,u)}\nonumber
\end{gather*}
and as usual
\[
y_{12}(q)=y(q,u_1-u_2) \quad \textrm{etc}.
\]
\begin{remark}Furthermore, the function $\omega(q,u)$ is actually independent of $q$, a property which will prove important later on,
and takes the following simple form:
\begin{gather}\label{omeg}
\omega(u)=\frac{\tf(u+\eta)\tf(u-2\eta)}{\tf(u-\eta)\tf(u+2\eta)}.
\end{gather}
\end{remark}

This identity can be proved by looking at transformation properties under $u\rightarrow u+1$, $u\rightarrow u+\tau$
of both sides of \eqref{omeg}.

\begin{remark} Notice also that $\omega(u)\omega(-u)=1$.
\end{remark}

Now we turn to the construction of the Bethe state. In the application of algebraic Bethe ansatz
to the $E_{\tau,\eta}(sl_2)$ elliptic quantum group the algebra contains a generator
(usually also denoted by $B(u)$) which acts as a creation operator. It also
enjoys the property $B(u)B(v)=B(v)B(u)$. This allows for the straightforward construction of the creation
operators $\Phi_n$ as
\[
\Phi_n(u_1,\ldots,u_n)=B(u_1)B(u_2)\cdots B(u_n),
\]
since we immediately have the property
\[
\Phi_n(u_1,\ldots,u_n)=\Phi_n(u_1,\ldots,u_{i-1},u_{i+1},u_i,u_{i+2},\ldots,u_n), \qquad i=1,2,\ldots,n-1.
\]

As it turns out, in the $E_{\tau,\eta}(so_3)$ case the creation operators are not
simple functions of the Lax matrix entries but they are complicated
functions of three generators $A_1(u)$, $B_1(u)$, $B_2(u)$ in the elliptic operator algebra.
This situation is analogous to that of the Izergin--Korepin model as described by Tarasov in~\cite{Ta}.

We give the following def\/inition for the creation operator.
\begin{definition}
Let $\Phi_n$ be def\/ined be the recurrence relation for $n\geq 0$:
\begin{gather*}
\Phi_n(u_1,\ldots,u_n)=B_1(u_1)\Phi_{n-1}(u_2,\ldots, u_n)\\
\phantom{\Phi_n(u_1,\ldots,u_n)=}{} -\sum_{j=2}^n\frac{\prod\limits_{k=2}^{j-1}\omega_{jk}}{y_{1j}(q)}
\prod_{\substack{k=2\\k\neq j}}^n z_{kj}(q+2\eta)\ B_2(u_1) \Phi_{n-2}(u_2,\ldots,\widehat{u_j},\ldots,u_n)A_1(u_j),
\end{gather*}
where $\Phi_0=1$, $\Phi_1(u_1)=B_1(u_1)$ and the hat means that that parameter is omitted.
\end{definition}

It may be useful to give explicitly the f\/irst three creation operators
\begin{gather*}
\Phi_1(u_1)=B_1(u_1),\\
\Phi_2(u_1,u_2)=B_1(u_1)B_1(u_2)-\frac{1}{y_{12}(q)}B_2(u_1)A_1(u_2),\\
\Phi_3(u_1,u_2,u_3)=B_1(u_1)B_1(u_2)B_1(u_3)-\frac{1}{y_{23}(q)}B_1(u_1)B_2(u_2)A_1(u_3)\\
\phantom{\Phi_3(u_1,u_2,u_3)=}{}-\frac{z_{32}(q+2\eta)}{y_{12}(q)}B_2(u_1)B_1(u_3)A_1(u_2)-\frac{\omega_{32}z_{23}(q+2\eta)}{y_{13}(q)}B_2(u_1)
B_1(u_2)A_1(u_3).
\end{gather*}

The Bethe vector is then not completely symmetric under the interchange of two neighboring spectral parameters
but verif\/ies the following property instead:
\begin{gather*}
\Phi_2(u_1,u_2)=\omega_{21}\Phi_2(u_2,u_1),\\
\Phi_3(u_1,u_2,u_3)=\omega_{21}\Phi_3(u_2,u_1,u_3)=\omega_{32}\Phi_3(u_1,u_3,u_2).
\end{gather*}

For general $n$ we prove the following theorem.
\begin{theorem}
$\Phi_n$ verif\/ies the following symmetry property:
\begin{gather}\label{symm}
\Phi_n(u_1,\ldots,u_n)=\omega_{i+1,i}\Phi_n(u_1,\ldots,u_{i-1},u_{i+1},u_i,u_{i+2},\ldots,u_n),
\qquad i=1,2,\ldots,n-1.\!\!\!
\end{gather}
\end{theorem}

\begin{proof}
The proof is by induction on $n$. The symmetry property is immediately proved for $i\neq 1$.
To verify for $i=1$, we have
to expand $\Phi_n$ by one more induction step:
\begin{gather*}
\Phi_n(u_1,\ldots,u_n)=B_1(u_1)B_1(u_2)\Phi_{n-2}(u_3,\ldots,u_n)\\
\qquad{}-\frac{\prod\limits_{k=3}^n z_{k2}(q+2\eta)}{y_{12}(q)}
B_2(u_1)\Phi_{n-2}(u_3,\ldots,u_n)A_1(u_2) \nonumber\\
\qquad{}-\sum_{j=3}^n \frac{\prod\limits_{k=3}^{j-1}\omega_{jk}}{y_{2j}(q)}
\prod_{\substack{k=3\\k\neq j}}^n z_{kj}(q+2\eta)
B_1(u_1)B_2(u_2)\Phi_{n-3}(u_3,\ldots,\widehat{u_j},\ldots,u_n)A_1(u_j) \nonumber\\
\qquad{}-\sum_{j=3}^n \frac{\omega_{j2}z_{2j}(q)\prod\limits_{k=3}^{j-1}
\omega_{jk}}{y_{1j}(q)}\prod_{\substack{k=3\\k\neq j}}^n
z_{kj}(q+2\eta) B_2(u_1)B_1(u_2)\Phi_{n-3}(u_3,\ldots,\widehat{u_j},\ldots,u_n) A_1(u_j) \nonumber\\
\qquad{}+\sum_{3\leq l \leq j \leq n}\left[ \frac{\omega_{j2}z_{2j}(q+2\eta)z_{lj}(q+2\eta)}{y_{1j}(q)y_{2l}(q+2\eta)}
+\omega_{lj}\frac{\omega_{l2}z_{2l}(q+2\eta)z_{jl}(q+2\eta)}{y_{1l}(q)y_{2j}(q+2\eta)} \right]\nonumber \\
\qquad{}\times \prod_{k=3}^{j-1}
\omega_{jk}\prod_{k=3}^{l-1}\omega_{lk}
 \prod_{\substack{k=3\\k \neq j, l}}^n z_{kj}(q+4\eta)z_{kl}(q+4\eta)
\frac{\tf(q+\eta)^2}{\tf(q-\eta)\tf(q+3\eta)} \nonumber\\ 
\qquad{}\times B_2(u_1)B_2(u_2)
\Phi_{n-4}(u_3,\ldots,\widehat{u_l},\ldots, \widehat{u_j},\ldots,u_n)A_1(u_l)A_1(u_j) \nonumber
\end{gather*}
then substitute into (\ref{symm})
and bring the right-hand-side to normal order of the
spectral parameters by using relations \eqref{crB1B1}--\eqref{crB2B1}. We f\/ind then that property (\ref{symm}) is
fulf\/illed provided the following identities hold true:
\begin{gather*}%\label{A1.4}
-\frac{\omega_{12}g_{21}}{y_{23}(q)\be_{21}(\eta,-q)}+\frac{\al_{21}(\eta,-q)}{\be_{21}(\eta,-q)y_{13}(q)}=
-\frac{\omega_{31}z_{13}(q+2\eta)}{y_{23}(q)}-\frac{\al_{31}(\eta,q+2\eta)}{\be_{31}(q+2\eta,\eta)y_{21}(q)}
\end{gather*}
and
\begin{gather*}
\omega_{12}\left(\frac{\omega_{42}z_{24}(q+2\eta)z_{34}(q+2\eta)}{y_{14}(q)y_{23}(q+2\eta)}+
\omega_{34}\frac{\omega_{32}z_{23}(q+2\eta)z_{43}(q+2\eta)}{y_{13}(q)y_{24}(q+2\eta)} \right) \\
\qquad{}-\left( \frac{\om_{41}z_{14}(q+2\eta)z_{34}(q+2\eta)}{y_{24}(q)y_{13}(q+2\eta)}+\frac{\om_{34}\om_{31}
z_{13}(q+2\eta)z_{43}(q+2\eta)}{y_{23}(q)y_{14}(q+2\eta)}\right) \\
\qquad{}+\frac{\om_{12}}{y_{12}(q)}\left( \frac{\de_{42}(-q-2\eta)}{\ga_{42}(q+2\eta,-q-2\eta)y_{43}(q)}+
\frac{z_{42}(q+2\eta)\al_{32}(\eta,q+2\eta)\om_{24}}{\be_{32}(q+2\eta,\eta)y_{24}(q+2\eta)} \right) \\
\qquad{}-\frac{1}{y_{21}(q)}\left(\frac{\de_{41}(-q-2\eta)}{\ga_{41}(q+2\eta,-q-2\eta)y_{43}(q)}+\frac{z_{41}(q+2\eta)
\al_{31}(\eta,q+2\eta)\om_{14}}{\be_{31}(q+2\eta,\eta)y_{14}(q+2\eta)}
\right)=0.
\end{gather*}

These identities can be verif\/ied by tedious calculations using once again quasiperiodicty pro\-per\-ties of
the $\tf$-function.
\end{proof}

\begin{remark} $\Phi_n$ contains the more familiar string of $B_1(u_1)\cdots B_1(u_n)$ with coef\/f\/icient $1$.
\end{remark}

It is straightforward to check the following relations using the commutation relations \eqref{RLLti}
\begin{gather*}
t(u)\Phi_1(u_1)= z_{1u}(q)B_1(u_1)A_1(u)-\frac{\al_{1u}(\eta,q)}{\be_{1u}(q,\eta)}B_1(u)A_1(u_1)\\
\qquad{}+\frac{z_{u1}(q)}{\omega_{u1}}B_1(u_1)A_2(u)-\frac{\al_{u1}(q,\eta)}{\be_{u1}(q,\eta)}B_1(u)A_2(u_1)
+\frac{1}{y_{u1}(q)}
B_3(u)A_1(u_1)\\
\qquad{}+\frac{\be_{u1}(\eta,-q)}{\ga_{u1}(q,-q)}B_1(u_1)A_3(u)-\frac{\ga_{u1}(q,\eta)}{\ga_{u1}(q,-q)}B_3(u)A_2(u_1),
\end{gather*}
for $n=2$
\begin{gather*}
t(u)\Phi_2(u_1,u_2)= z_{1u}(q)z_{2u}(q)\Phi_2(u_1,u_2)A_1(u)+\frac{z_{u1}(q)z_{u2}(q-2\eta)}{\omega_{u1}\omega_{u2}}
\Phi_2(u_1,u_2)A_2(u)\\
{}+\frac{\be_{u1}(\eta,-q)\be_{u2}(\eta,-q)}{\ga_{u1}(q-2\eta,-q+2\eta)\ga_{u2}(q-2\eta,-q+2\eta)}
\Phi_2(u_1,u_2)A_3(u)\\
{}+\left(-\frac{z_{1u}(q)\al_{2u}(\eta,q)}{\be_{2u}(q,\eta)}+\frac{\al_{1u}(\eta,q)\al_{21}(\eta,q)\omega_{1u}}
{\be_{1u}(q,\eta)\be_{21}(q,\eta)}
-\frac{\ga_{1u}(\eta,-q)}{y_{12}(q-2\eta)ga_{1u}(q,-q)}\right)\Phi_2(u_1,u)A_1(u_2)\\
{}+\left(-\frac{z_{u1}(q)\al_{u2}(q-2\eta,\eta)}{\omega_{u1}\be_{u2}(q-2\eta,\eta)}+
\frac{\al_{u1}(q,\eta)\al_{12}(q-2\eta,\eta)}{\be_{u1}(q,\eta)\be_{12}(q-2\eta,\eta)\omega_{u1}} \right)
\Phi_2(u_1,u)A_2(u_2)\\
{}-\frac{z_{21}(q)\al_{1u}(\eta,q)}{\be_{1u}(q,\eta)}\Phi_2(u,u_2)A_1(u_1)- \frac{\al_{u1}(q,\eta)z_{12}(q-2\eta)}
{\be_{u1}(q,\eta)\omega_{12}}\Phi_2(u,u_2)A_2(u_1)\\
+\left(-\frac{z_{21}(q)\al_{1u}(\eta,q)}{y_{u2}(q)\be_{1u}(q,\eta)}+\frac{\al_{1u}(\eta,q)\al_{21}(\eta,q)}
{\be_{1u}(q,\eta)\be_{1u}(q,\eta)y_{u1}(q)}
-\frac{\ga_{1u}(\eta,-q)}{\ga_{1u}(q,-q)y_{12}(q-2\eta)y_{u1}(q)}\right.\\
\left.{}+\frac{\de_{1u}(-q)}{\ga_{1u}(q,-q)y_{12}(q-2\eta)}
\right)B_2(u)A_1(u_1)A_1(u_2)\\
\times\left( \frac{z_{u1}(q)\al_{u1}(q,\eta)}{\omega_{u1}y_{u2}(q)\be_{u1}(\eta,-q)}-\frac{\al_{u1}(q,\eta)g_{12}}
{\be_{u1}(q,\eta)\omega_{12}y_{u2}(q)\be_{12}(q+2\eta,\eta)}+\frac{\al_{u1}(q,\eta)\al_{12}(q+2\eta,\eta)}
{\be_{u1}(q,\eta)y_{u1}(q)\be_{12}(q+2\eta,\eta)}\right.\!\!\\
\left.{}- \frac{\al_{u1}(q,\eta)\al_{u1}(\eta,-q)}{y_{12}(q)\be_{u1}(\eta,-q)
\be_{u1}(q,\eta)} \right)B_2(u)A_2(u_1)A_1(u_2)\\
 {}\times \frac{1}{\ga_{u1}(q,-q)}\left(\frac{\de_{u1}(q)}{y_{12}(q-2\eta,-q+2\eta)}-\frac{\al_{u1}(q,\eta)}{y_{u2}(q-2\eta)}
\right)B_2(u)A_2(u_1)A_2(u_2)\\
{}\times\left( \frac{g_{u1}}{\be_{u1}(\eta,-q)\omega_{u1}y_{u2}(q)}-\frac{\al_{21}(\eta,q+2\eta)}{y_{u1}(q)
\be_{21}(q+2\eta,\eta)}
-\frac{\al_{u1}(\eta,-q)}{y_{12}(q)\be_{u1}(\eta,-q)}\right)B_3(u)B_1(u_1)A_1(u_2)\\
{}+\frac{\al_{12}(q,\eta)}{y_{u1}(q,-q)\be_{12}(q,\eta)}B_3(u)B_1(u_1)A_2(u_2)\\
\frac{z_{21}(q+2\eta)}{y_{u1}(q)}B_3(u)B_1(u_1)A_2(u_1)-\frac{z_{12}(q)}{y_{u1}(q)\omega_{12}}B_3(u)B_1(u_2)A_2(u_1)\\
{}+\textrm{terms ending with}\  C.
\end{gather*}

The cancelation of the unwanted terms is ensured for $n=2$ if $u_1$ and $u_2$ are solutions of the following
Bethe equations.
\begin{gather*}
\frac{a_1(u_i)}{a_2(u_i,q)}=\prod_{\substack{j=1\\j\neq i}}^{2}\frac{\tf(u_{ij}-\eta)}{\tf(u_{ij}+\eta)}\times
\frac{\tf(q-3\eta)^2 }{\tf(q-\eta)\tf(q-5\eta)}\times \frac{f(q)}{f(q-2\eta)}, \qquad i=1,2.
\end{gather*}

The role of the function $f(q)$ becomes clear in this context. It has to be chosen so as to eliminate
the $q$-dependence from the Bethe equation.

These results suggest that $\Phi_n$ for all $n$ are the correct choice of creation operators of the
corresponding Bethe states. The complete proof of the general case will be given elsewhere.

\section{Conclusion}

In this paper we def\/ined the elliptic quantum group $E_{\tau,\eta}(so_3)$ along the lines described
by Felder in \cite{Fe}. Although dynamical, the $R$-matrix appearing in the exchange relations
has a matrix form similar to that of the Izergin--Korepin model. Lax operator, operator algebra and families
of commuting transfer matrices are def\/ined in complete analogy with the $E_{\tau,\eta}(sl_2)$ case.

Our aim was to apply algebraic Bethe ansatz method in this setting. We have obtained a~recurrence
relation for the creation operators and have proved that these operators have a~certain symmetry property
under the interchange of two adjacent spectral parameters. Both the form of the recurrence relation and
this symmetry property are an elliptic generalization of Tarasov's results in~\cite{Ta}. Finally,
we have obtained the Bethe equations and the eigenvalues for the 2-magnon state.
The Bethe
equations for the general n-magnon state
will be published later~\cite{MaNa2}.

\subsection*{Acknowledgements}

We wish to thank Petr %Petrovich 
Kulish for illuminating discussions.
We are also grateful to the organizers of the GEOMIS workshop in Coimbra for letting
us present these results. This work was supported by the project POCI/MAT/58452/2004,
in addition to that Z.~Nagy benef\/ited from
the FCT grant SFRH/BPD/25310/2005. The manuscript was f\/inished
in the hospitable environment of the Solvay Institute and ULB, Brussels.

\pdfbookmark[1]{References}{ref}

\LastPageEnding

\end{document}